\documentclass[preprint,12pt]{elsarticle}

\usepackage{amssymb, amsmath, amsthm}
\usepackage[colorlinks]{hyperref}

\journal{Stochastic Processes and their Applications}

\begin{document}

\frenchspacing
\renewcommand{\epsilon}{\varepsilon}
\renewcommand{\phi}{\varphi}

\newtheorem{thm}{Theorem}[section]
\newtheorem{prop}[thm]{Proposition}
\newtheorem{lem}[thm]{Lemma}
\newtheorem{dfn}[thm]{Definition}
\newtheorem{cor}[thm]{Corollary}

\theoremstyle{definition}
\newtheorem{ex}[thm]{Example}
\newtheorem{rem}[thm]{Remark}
\newtheorem{alg}[thm]{Algorithm}

\renewenvironment{proof}{{\bf Proof:}}
{\hspace*{\stretch{1}} $\square$}

\begin{frontmatter}

%% Title, authors and addresses

%% use the tnoteref command within \title for footnotes;
%% use the tnotetext command for the associated footnote;
%% use the fnref command within \author or \address for footnotes;
%% use the fntext command for the associated footnote;
%% use the corref command within \author for corresponding author footnotes;
%% use the cortext command for the associated footnote;
%% use the ead command for the email address,
%% and the form \ead[url] for the home page:
%%
%% \title{Title\tnoteref{label1}}
%% \tnotetext[label1]{}
%% \author{Name\corref{cor1}\fnref{label2}}
%% \ead{email address}
%% \ead[url]{home page}
%% \fntext[label2]{}
%% \cortext[cor1]{}
%% \address{Address\fnref{label3}}
%% \fntext[label3]{}

\title{Penalized maximum likelihood estimation for generalized linear point processes}

%% use optional labels to link authors explicitly to addresses:
%% \author[label1,label2]{<author name>}
%% \address[label1]{<address>}
%% \address[label2]{<address>}

\author{Niels Richard Hansen}
\ead{Niels.R.Hansen@math.ku.dk}

\address{Department of Mathematical Sciences, University
of Copenhagen, Universitetsparken 5, 2100 Copenhagen \O,
Denmark.}

\begin{abstract}
  A generalized linear point process is specified in terms of an
  intensity that depends upon a linear predictor process through a
  fixed non-linear function. We present a framework where the linear
  predictor is parametrized by a Banach space and give results on
  G\^ateaux differentiability of the log-likelihood. Of
  particular interest is when the intensity is expressed in terms of a
  linear filter parametrized by a Sobolev space. Using that the
  Sobolev spaces are reproducing kernel Hilbert spaces we derive
  results on the representation of the penalized maximum likelihood
  estimator in a special case and the gradient of the negative
  log-likelihood in general. The latter is used to develop a descent
  algorithm in the Sobolev space. We conclude the paper by extensions
  to multivariate and additive model specifications. The methods are
  implemented in the R-package \verb+ppstat+.
\end{abstract}

\begin{keyword}
%% keywords here, in the form: keyword \sep keyword

%% MSC codes here, in the form: \MSC code \sep code
%% or \MSC[2008] code \sep code (2000 is the default)

\end{keyword}

\end{frontmatter}

%%
%% Start line numbering here if you want
%%
% \linenumbers

\section{Introduction}

In this paper we aim at combining likelihood based inference for
stochastic processes with non-parametric regression methods. In
particular, we discuss estimation of smooth functional components in linear filters
that enter in the specification of a point process model. The results
were inspired by applications of multivariate point process models to the
modeling of the occurrences of transcription regulatory elements along
the genome and the activity of collections of neurons. 

There are many important applications of one-dimensional point process
models such as models of queuing and telecommunication systems,
\citep{Asmussen:2003}, insurance claims, \citep{Mikosch:2004},
earthquakes, \citep{Ogata:1986}, \citep{Ogata:2003}, neuronal
activity, \citep{Brillinger:1992}, \citep{Paninski:2004},
high-frequency financial activity, \citep{Hautsch:2004}, and
occurrences of DNA motifs, \citep{Gusto:2005}, \citep{Reynard:2009},
just to mention some. Andersen et al.,
\citep{AndersenBorganGillKeiding:1993}, give a general treatment of
statistics for point process model -- with a focus on 
applications in event history analysis. See also
\citep{FlemingHarrington:1991} or \citep{Karr:1991} for general
introductions to statistics for point processes. Some recent
applications of multivariate point processes, a.k.a.  marked point
processes, include our integrated analysis of ChIP-seq data,
\citep{Carstensen:2010}, the modeling of multivariate neuron spike
data, \citep{Pillow:2008}, \citep{Masud:2011}, and stochastic kinetic
modeling, \citep{Bowsher:2010}.

  In our work on genomic organization of transcription regulatory
  elements based on ChIP-chip and ChIP-seq data,
  \citep{Carstensen:2010}, we were inspired by the use of linear
  Hawkes processes in \citep{Gusto:2005}, and the general class of
  multivariate, non-linear Hawkes processes, as treated in
  \citep{Bremaud:1996}. We developed a first version of the R-package
  \verb+ppstat+ for the likelihood based analysis using non-linear
  Hawkes processes. The Hawkes models share a structural similarity
  with generalized linear models, and it is possible to carry out the
  practical computations using Poisson regression methods. The
  terminology of a generalized linear point process model has,
  furthermore, been used recently for various Hawkes-like models of
  spike trains for neurons, \citep{Paninski:2004},
  \citep{Pillow:2008}, \citep{Toyoizumi:2009}. The models considered
  in \citep{Pillow:2008} for multivariate spike trains share many
  similarities with our models of the occurrences of multiple
  transcription regulatory elements. In particular, the use of basis
  expansions for estimation of functional components, which may be
  combined with regularization in terms of penalized
  maximum-likelihood estimation. In \citep{Pillow:2008} the basis
  functions chosen were raised cosines with a log-time transformation,
  whereas we used B-splines in \citep{Carstensen:2010}. 

  We found it useful to give a general definition of a
  \emph{generalized linear point process model} as a process where the
  intensity is linked to a predictor process, which is linear
  in the unknown parameters, and where this linear predictor process
  potentially depends on the internal history of the point process as
  well as additional covariate processes. The R-package \verb+ppstat+
  has been developed for likelihood based analysis of data
  from multivariate point processes. The package handles, in
  particular, the non-linear Hawkes processes where intensities are
  given in terms of a non-linear function of linear filters with
  filter functions given via basis expansions. Its usage is documented
  in detail elsewhere, see
  \href{http://www.math.ku.dk/~richard/ppstat/}{\url{http://www.math.ku.dk/~richard/ppstat/}}.
  See also \cite{Hansen:2013b} for computational details.

  The focus of the present paper is on the theoretical framework for
  the computation of penalized maximum-likelihood estimators of
  functional parameters in a one-dimensional point process setup. For
  a treatment of sampling properties of penalized maximum-likelihood
  estimators see \citep{Cox:1990}. We show how a particular set of
  basis functions appears as the solution of a more abstractly
  formulated problem. We have the classical result on smoothing
  splines in mind, which says that the solution of a
  roughness-penalized least squares problem is a spline, see Theorem
  2.4 in \citep{Green:1994}. We first introduce the framework of
  generalized linear point process models parametrized by a Banach
  space, and we give general results on derivatives of the 
  log-likelihood function. Then we restrict attention to a particular
  class of linear filters parametrized by Sobolev spaces that includes
  the non-linear Hawkes processes as a special case.  We show two main
  results for this class of models. The first result we show is
  similar to the result on smoothing splines, and it states that the
  penalized maximum-likelihood estimator in a special case is found in
  a finite-dimensional space spanned by an explicit set of basis
  functions. For the linear Hawkes process the solution is a
  spline. The second result is different. For the general model class
  considered we do not find an explicit finite-dimensional
  basis. Instead we derive an infinite-dimensional gradient, which
  suggests an iterative algorithm, and we establish a convergence
  result for this algorithm. The algorithm can be interpreted as a
  sequence of finite-dimensional subspace approximations. We exploit 
  that Sobolev spaces are reproducing kernel
  Hilbert spaces, and that the likelihood in the special case and the
  gradient of the log-likelihood in general are given in
  terms of continuous linear functionals. These functionals are
  expressed as stochastic integrals of integrands from a Sobolev space. In a
  regression context the linear functionals considered are typically
  simple point evaluations, which are trivially continuous. In the
  context of the present paper it is more involved to establish
  continuity, and we use specific properties of Sobolev spaces as well as
  their general properties as reproducing kernel Hilbert spaces.

\section{Setup}

We let $(\Omega,\mathcal{F},(\mathcal{F}_t)_{t \geq 0},P)$ be a filtered probability space -- a stochastic basis -- 
where the filtration is
assumed to be right continuous. We will, in addition, assume that
$(N_t)_{t \geq 0}$ is an adapted counting process, which, under $P$,
is a homogeneous Poisson process with rate 1. 

If $(\lambda_t)_{t \geq 0}$ is a positive, predictable process we 
define the \emph{likelihood process}
\begin{equation} \label{eq:likelihoodprocess}
\mathcal{L}_t = \exp\left(t + \int_{0}^t \log \lambda_s
\mathrm{d}N_s -\Lambda_t\right), \qquad \Lambda_t = \int_{0}^t \lambda_s \mathrm{d}s.
\end{equation}
We will assume that $\Lambda_t < \infty$ $P$-a.s., in which case
$(\mathcal{L}_t)_{t \geq 0}$ is a $P$-local martingale and
a $P$-supermartingale with $\mathbb{E}_P(\mathcal{L}_t) \leq 1$ for
all $t \geq 0$, see Theorem VI.T2, \citep{Bremaud:1981}. If
$\mathbb{E}_P(\mathcal{L}_t) = 1$ we can define a probability measure
$Q_t$ on $\mathcal{F}$ by taking $\mathcal{L}_t$ to be the
Radon-Nikodym derivative of $Q_t$ w.r.t. $P$.  That is,
\begin{equation} \label{eq:QtDef1}
Q_t = \mathcal{L}_t \cdot P.
\end{equation}
We note that $\mathbb{E}_P(\mathcal{L}_t) = 1$ if and only
if $(\mathcal{L}_s)_{0 \leq s \leq t}$ is a true
$P$-martingale. If $\mathbb{E}_P(\mathcal{L}_t) < 1$ we cannot define
a probability measure $Q_t$ on the abstract space $(\Omega,\mathcal{F})$ by
(\ref{eq:QtDef1}). With a canonical choice of $\Omega$
it is always possible to construct a measure $Q_t$ such that 
\begin{equation}  \label{eq:QtDef}
Q_t = \mathcal{L}_t \cdot P + Q_t^{\perp}
\end{equation}
where $Q_t^{\perp}(N_t < \infty) = 0$, see \citep{Jacod:1975} or
Theorem 5.2.1(ii), \citep{Jacobsen:2006}. General conditions assuring
that $\mathbb{E}_P(\mathcal{L}_t) = 1$ can be found in \cite{Sokol:2012}. Though
it is important to be able to
decide if the likelihood process is a true martingale, it plays no
role for the results and computations in the present paper. 

Throughout we will fix an observation window $[0,t]$ and assume that
we have observed a non-exploding realization of $(N_s)_{0 \leq s \leq t}$ under a $Q_t$-measure fulfilling
(\ref{eq:QtDef}). The process
$(\lambda_s)_{0 \leq s \leq t}$ is called the (predictable) intensity
process for the counting process $(N_s)_{0 \leq s \leq t}$ under
$Q_t$. The integrated intensity, $(\Lambda_s)_{0 \leq s \leq t}$, is
the compensator, and if $\mathbb{E}_P(\mathcal{L}_t) =1$ 
the process $M_s = N_s - \Lambda_s$ for $s \in [0,t]$ is a
$Q_t$-martingale, see Theorem VI.T3, \citep{Bremaud:1981}.

We will study models where the intensity is parametrized by a Banach
space valued parameter. Let $V$ denote a Banach space
with $V^*$ its dual space of continuous
linear functionals. We equip 
$V^*$ with the $\sigma$-algebra\footnote{If $V$ is separable the dual space $V^*$ is separable and second
countable in the weak$^*$-topology in which case the $\sigma$-algebra coincides with the weak$^*$ Borel $\sigma$-algebra.} generated by the linear functionals 
$$x \mapsto x\beta$$
for $\beta \in V$. We observe that if $X(\omega)$ is a
linear functional on $V$ it belongs to $V^*$ if and only if $\beta
\mapsto X(\omega)\beta$ is continuous, and if $X(\omega) \in V^*$ for
all $\omega$ then $X$ is measurable as a
map $X: \Omega \to V^*$ if and only if $\omega \mapsto X(\omega)\beta$
is measureable for all $\beta \in V$. A stochastic
process $(X_s)_{0 \leq s\leq t}$ with values in $V^*$ is thus 
\emph{adapted} if and only if $(X_s\beta)_{0 \leq s\leq t}$ is
adapted. 

We say that a stochastic process $(X_s)_{0 \leq s\leq t}$ with values in
$V^*$ is continuous from the left (right) and has limits from the right
(left) if this holds for $(X_s\beta)_{0 \leq s\leq t}$ for all $\beta
\in V$. Thus these continuity properties of $s \mapsto X_s$ from
$[0,t]$ into $V^*$ are with respect to the weak$^*$-topology on $V^*$.

\begin{dfn} \label{dfn:ppglm} Let $(X_s)_{0 \leq s\leq t}$ be an
  adapted process with values in $V^*$, continuous from the left and with
  right limits. Let $\phi : D \to
  [0,\infty)$ for $D \subseteq \mathbb{R}$ be continuous and let 
$$\Theta(D) = \{ \beta \in V \mid X_{s}\beta \in D \textrm{
  for all } s \in [0,t] \}.$$
  A generalized linear point
  process model on $[0,t]$ is a point process on $[0,t]$ parametrized
  by $\Theta(D)$ such that for $\beta \in \Theta(D)$ the point
  process has intensity 
$$\lambda_s = \phi(X_{s}\beta)$$
for $s \in [0,t]$.
\end{dfn}

Continuity from the left and adaptedness ensures predictability of the
intensity, cf. Definition 2.1 in \citep{Jacod:2003}. Requiring finite limits
from the right ensures boundedness ($\omega$-wise) of $s \mapsto
\phi(X_{s}\beta)$ on $[0,t]$
and thus that $\int_0^t \phi(X_{s}\beta) \mathrm{d} s < \infty$. 

We call $(X_{s} \beta )_{0 \leq s \leq t}$ the \emph{linear predictor
  process}, which can be be interpreted as a predictable filter of the
Banach space valued process $(X_s)_{0 \leq s \leq t}$. The possible
filters are parametrized by $\beta \in \Theta(D)$, and the objective,
from a statistical point of view, is the estimation of $\beta$. The
definition includes the possibility of $V = C_b(\mathbb{R})$, the
space of bounded continuous functions equipped with the uniform norm,
and $X_{s}\beta = \beta(X_{s})$ for a real valued predictable process
$X$. This evaluation filter is a non-linear filter in $X_{s}$ but
linear in $\beta$. A particular example is the inhomogeneous
Poisson process obtained by taking $X_s = s$.

Our main focus, as presented in Section \ref{sec:linearfilter}, is to
the case where $V$ is a reproducing kernel Hilbert space, and where $X_s$ is given in terms of stochastic integration w.r.t. an
ordinary real valued stochastic process. These filters will be linear
filters in the stochastic process. 

Note that if $\phi$ is one-to-one
with inverse $m = \phi^{-1} : \phi(D) \to D$ then 
 $$X_{s}\beta = m(\lambda_s).$$
 Drawing an analogy to ordinary generalized linear models it seems
 natural at this point to call $m$ the link function -- it transforms
 the intensity process into a process that is linear in the parameter
 $\beta$. With this terminology we would call $\phi$ the inverse link
 function. However, there is no reason to
 require $\phi$ to be one-to-one in general, and we will not use the terminology.  

When the likelihood process is a martingale it is evident from (\ref{eq:likelihoodprocess}) that as a statistical model with parameter space
$\Theta(D) \subseteq V$ the negative log-likelihood function
for observing $(N_s)_{0 \leq s \leq t}$ is 
\begin{equation}
  \ell_t(\beta) = \int_0^t \phi(X_{s} \beta) \mathrm{d}s
-  \int_0^t  \log(\phi(X_{s} \beta)) \mathrm{d}N_s
\end{equation}
for $\beta \in \Theta(D)$. Strictly speaking, $\ell_t$ is only a true
negative log-likelihood if $\mathbb{E}_P(\mathcal{L}_t) = 1$, but for
non-exploding data $\ell_t$ actually encodes all pairwise likelihood
comparisons even if the measures are not equivalent. Anyway, the
concerns of the present paper are representations and computations of
the penalized maximum-likelihood estimator based on $\ell_t$, in
which case it plays no role whether $\ell_t$ is a true negative log-likelihood.
As a final remark we note that the negative log-likelihood
function is convex as a function of $\beta$ if $\phi$ is convex and
log-concave. 

%  and we will focus on the penalty functions of the form 
% $J(\beta) = \lambda r(||\beta||)$ where $\lambda > 0$ is a tuning
% parameter and $r$ is an increasing function. If $V = \mathbb{R}^{d'}$ 
% and $B_{s,1},\ldots, B_{s,d'}$ are $C^2$ basis
% functions in $s$ we identify $\beta \in  \mathbb{R}^{d'}$ with the function 
% $$\beta(s) = \sum_{i=1}^{d'} \beta_i B_{s, i}$$ 
% and an important penalty function on $ \mathbb{R}^{d'}$
% is given by the quadratic form $K$ where 
% $$K_{i,j} = \int_0^t \partial_s^2 B_{s,i} \partial_s^2 B_{s,j} \mathrm{d} s,$$
% in which case
% \begin{eqnarray*} 
% J(\beta)  = \lambda \beta_{\alpha}^T K \beta_{\alpha} = \lambda \int_0^t
% \left( \partial_s^2 \beta(s) \right)^2 \mathrm{d} s.
% \end{eqnarray*}

% The generalized linear point process models have value when
% $V$ is finite dimensional, but we emphasize that models of 
% greater generality fit into the model
% class above for a suitable choice of $X$-process. This is at least true 
% from a practical point of view where finite basis expansions can be
% used to approximate non-parametric components, and we also show one
% result in Section \ref{sec:linearfilter} where 
% penalized maximum-likelihood estimation in an infinite dimensional
% function space reduces 
%  to penalized maximum likelihood estimation for a generalized linear
%  point process model with a finite dimensional parameter space. 

\section{Results}

Before turning to more concrete models we give general
results on differentiation of the
negative log-likelihood. First under the assumption that $\phi$ is
suitably differentiable, but subsequently illustrating that the
time-integral in the negative log-likelihood can smooth out
non-differentiabilities in $\phi$. All proofs are postponed to Section \ref{sec:proofs}.

\begin{prop} \label{prop:diff} If $\phi$ is $C^1$ on $D$, and
  $l_t(\beta) < \infty$ for $\beta \in \Theta(D)^{\circ}$ then $l_t$
  is G\^ateaux differentiable in $\beta$ with 
\begin{equation}
  D\ell_t(\beta) = \int_{0}^t \phi'(X_s \beta) X_s \mathrm{d}s
-  \int_{0}^t \frac{\phi'(X_s \beta)}{\phi(X_s \beta)}  X_s \mathrm{d}N_s.
\end{equation}
Moreover, if $\phi$ is $C^2$ the second G\^ateaux derivative is 
\begin{eqnarray} \nonumber
 D^2\ell_t(\beta) &= &\int_0^t \phi''(X_s \beta) X_s
 \otimes X_s \mathrm{d}s \\ && 
-  \int_0^t \frac{\phi''(X_s \beta) \phi(X_s \beta) -
  \phi'(X^T_{s-} \beta)^2}{\phi(X_s \beta)^2}  X_s \otimes X_s \mathrm{d}N_s.
\end{eqnarray}
\end{prop}

The integrals above are to be interpreted as weak, or Pettis,
integrals. From the formulas it follows that $D\ell_t(\beta)$ is
linear and $D^2\ell_t(\beta)$ is bilinear. However, without further
assumptions on $X_s$ neither needs to be continuous. Continuity
follows if $||X_s||$ can be bounded ($\omega$-wise) as a function of
$s$. This is one of the main questions we deal with in the context of
Section \ref{sec:linearfilter} -- specifically we derive a gradient
representation of the derivative in a reproducing kernel Hilbert space
by proving continuity of $D\ell_t(\beta)$. Knowledge of the second
derivative is used for quadratic approximations of the negative
log-likelihood, in particular in relation to Algorithm \ref{alg:1} and
the iterative optimization over finite dimensional subspaces.

Simple formulas are obtained with $\phi(x) = x$, but this choice of
$\phi$ puts an often inconvenient restriction on the parameter space
to ensure that the intensity stays positive. This can be circumvented
by taking $\phi(x) = x_+$, but then the formulas above
break down -- in particular for the second derivative.  A possible
workaround is to modify $\phi$ locally around $0$ to make it twice
continuously differentiable. It is, however, not obvious that the
resulting formulas for the derivative are numerically stable and play
together with the time discretization that eventually must be used for
computing the time integral. We show that if
$s \mapsto X_s\beta$ has a finite number of roots and is locally smooth
around the roots then the time integral smoothes the negative log-likelihood to make
it twice differentiable. 

\begin{prop} \label{prop:diff2} Take $\phi(x) = x_+$ and assume that $\beta \in
  \Theta(D)^{\circ}$ is such that $l_t(\beta) < \infty$ and
  $s \mapsto X_s\beta$ has a finite number of roots in
  $[0,t]$. Then $\ell_t$ is G\^ateaux differentiable with 
\begin{equation} \nonumber
  D\ell_t(\beta) = \int_{0}^t 1(X_s \beta > 0) X_s \mathrm{d}s
-  \int_{0}^t \frac{1}{X_s \beta}  X_s \mathrm{d}N_s.
\end{equation}
Moreover, if there are neighborhoods of the roots $s_1, \ldots, s_n$
in which the sample paths of $X_{s}\beta$ and $X_{s}\gamma$ are
$C^1$, and $\partial_s
X_{s_i}\beta \neq 0$ for $i = 1, \ldots, n$, then the second
G\^ateaux derivative in $(\rho, \gamma)$ is 
\begin{eqnarray} \nonumber
 D^2\ell_t(\beta)(\rho, \gamma) &= & \sum_{i=1}^n 
 \frac{1}{|\partial_s X_{s_i}\beta|}  X_{s_i}\rho X_{s_i}\gamma 
+  \int_0^t \frac{1}{(X_s \beta)^2}  X_s\rho X_s\gamma \mathrm{d}N_s.
\end{eqnarray}
\end{prop}

\subsection{Linear filters from stochastic integration} \label{sec:linearfilter}

Let $g: [0,\infty) \to \mathbb{R}$ be a 
measurable function and $(Z_s)_{0 \leq s \leq t}$ a c\`adl\`ag process. If $g$
is locally bounded and $Z$ is a semi-martingale the stochastic process
$$\int_0^{s-} g(s-u) \mathrm{d}Z_{u}$$
is a well defined process. The process is sometimes called a homogeneous linear
filter or a moving average. 

We will need to interpret the stochastic integral above as a
stochastic process with values in a dual space. Since stochastic
integrals are usually not defined pathwisely, it is, in fact, not
obvious that
\begin{equation} \label{eq:Xdef}
g \mapsto X_sg := \int_0^{s-} g(s-u) \mathrm{d}Z_{u}
\end{equation}
for a fixed sample path is even a well defined linear functional --
let alone continuous. If we take the parameter space for $g$ to be $V
= W^{m,2}$, that is, $V$ is the Sobolev space of functions on $[0,t]$ that
are $m$ times weakly differentiable with the $m$'th derivative in
$L_2([0,t])$, then $g$ is weakly
differentiable with $L_2$-derivative for $m \geq 1$. Hence, for $Z$ a
semi-martingale, we have by integration by parts that 
\begin{eqnarray} \label{eq:intebyparts}
\int_0^{s-} h(u) \mathrm{d}Z_{u} & = & h(s)Z_{s-} - h(0)Z_{0}- \int_0^s Z_{u}  h'(u) \mathrm{d}u  
\end{eqnarray}
for $h \in W^{m,2}$. This equality is in general valid up to evanescence. The
right hand side is pathwisely well defined, and we use this as the
pathwise definition of the stochastic integral of $h \in W^{m,2}$
w.r.t. a c\`adl\`ag process $Z$. 
The integral then becomes a linear functional in $h$ for a concrete realization of
the $Z$-process, and by Corollary
\ref{cor:fundcont} $X_s$ is a continuous linear
functional. Thus $(X_s)_{0 \leq s \leq t}$ is a stochastic
process with values in $V^*$. Lemma \ref{lem:caglad} shows, moreover,
that $(X_s)_{0   \leq s \leq t}$ is continuous from the left with
limits from the right.

If the function $\phi: D \to [0,\infty)$ is
given we find that $\Theta(D)$ consists of those $g$ such that 
\begin{equation} \label{eq:pardomaincond}
 \int_0^{s-} g(s-u) \mathrm{d}Z_{u} \in D \
\textrm{for all } s \in [0,t].
\end{equation}
The particular case of interest with $D \neq \mathbb{R}$ is $D =
[0,\infty)$ and $Z$ an increasing process, e.g. a counting
process. In this case $g \in \Theta([0,\infty))$ if $g \geq 0$.

The Sobolev space $W^{m,2}$ can be equipped with several inner
products that give rise to equivalent norms and turn the space into a
reproducing kernel Hilbert space, \citep{Wahba:1990}, \citep{Berlinet:2004}. For each inner product there is an associated kernel, the
reproducing kernel, and we assume here that one inner product is
chosen with the corresponding norm denoted $||\cdot||$ and corresponding kernel denoted $R:[0,t]\times
[0,t] \to \mathbb{R}$.  Moreover, we fix $\gamma_1,\ldots,\gamma_l \in
W^{m,2}$ and denote by $P$ the orthogonal projection onto  
$\text{span}\{\gamma_1,\ldots,\gamma_l\}^{\perp}$. One of the defining
properties of the kernel $R$ is that for fixed $s \in [0,t]$,
$R(s,\cdot) \in W^{m,2}$, hence $PR(s,\cdot)$ is a well defined
function. This give rise to the projected kernel, which we denote $R^1
= PR$. The penalized negative log-likelihood function we consider is 
\begin{equation} \label{eq:minfunc}
l_t(g) + \lambda  ||Pg||^2 
\end{equation}
for $g \in \Theta(D)$ and $\lambda > 0$ where 
$$\ell_t(g) = \int_0^t \phi\left(  \int_0^{s-}
    g(s-u) \mathrm{d}Z_{u}  \right) \mathrm{d}s
-  \int_0^t \log( \phi\left( \int_0^{s-} g(s-u)
    \mathrm{d}Z_{u}\right)) \mathrm{d}N_s.$$
With $\tau_1, \ldots, \tau_{N_{t}}$ denoting the jump times for the
counting process $(N_s)_{0 \leq s \leq t}$ we can state one of the main theorems. 

\begin{thm} \label{thm:1} If $\phi(x) = x+d$ with domain
  $D=[-d,\infty)$ then a minimizer of (\ref{eq:minfunc}) over 
  $\Theta(D) \subseteq W^{m,2}$, $m \geq 1$, belongs to the finite dimensional subspace
  of $W^{m,2}$ spanned by the functions $\gamma_1,\ldots,\gamma_l$, the functions
  $$h_{i}(r) = \int_0^{\tau_i-} R^1(\tau_i-u,r) \mathrm{d}Z_{u}$$
for $i=1,\ldots,N_t$ together with the function 
$$f(r) = \int_0^t  \int_0^{s-} R^1(s-u,r)  \ \mathrm{d}Z_{u} \mathrm{d}s.$$
\end{thm}

\begin{rem} \label{rem:prac} A practical consequence of Theorem
  \ref{thm:1} is that the estimation problem reduces to a 
finite dimensional optimization problem over the space spanned by the
$l+1+N_t$ dimensional vector formed by combining $\gamma_1, \ldots, \gamma_l$, $f$ and $h_{i}$,  $i=1,\ldots,N_t$. For the concrete realization we 
may of course choose whichever basis that is most convenient for this
function space. For the practical computation of $f$ we note that by
Lemma \ref{lem:fubini} we can interchange the order of the integrations so
that
\begin{equation} \label{eq:fbasisinterchange}
f(r) = \int_0^t  \int_{u}^{t} R^1(s-u,r) \ \mathrm{d}s
\mathrm{d}Z_{u}.
\end{equation}
A detailed example is worked out in
Section \ref{sec:examples}.
\end{rem} 

\begin{rem} It is a common trick to construct a model conditionally on the entire 
outcome of a process $(Z_s)_{0 \leq s \leq t}$ by assuring that $Z_s$ is $\mathcal{F}_0$-measurable 
for all $s \in [0,t]$. In this case the process  
$$\int_0^t g(|s-u|) \mathrm{d}Z_{u}$$
for $s \in [0,t]$ becomes predictable. 
Theorem \ref{thm:1} still holds with the modification that 
$$h_{i}(r) = \int_0^t R^1(|\tau_i-u|,r) \mathrm{d}Z_{u}$$
for $i=1,\ldots,N_t$ and 
$$f(r) = \int_0^t  \int_0^t R^1(|s-u|,r) \mathrm{d}Z_{u} \mathrm{d}s.$$
When we model events that happen in time it is most natural that the
intensity at a given time $t$ only depends on the behavior of the
$Z$-process up to just before $t$. This corresponds to the formulation
chosen in Theorem \ref{thm:1}. However, if we model events in a
one-dimensional space it is often more natural to take the approach in
this remark.
\end{rem}  

If $\phi$ is not an affine function, we
cannot compute an explicit finite dimensional subspace. Instead, we
compute the gradient of the negative log-likelihood function. 

\begin{prop} \label{thm:2} If $\phi$ is continuously differentiable
  and $g \in \Theta(D)^{\circ}$ we define  
 $\eta_i$ for $i=1,\ldots,N_t$ as
$$\eta_i(r) = \int_0^{\tau_i-} R(\tau_i-u,r) \mathrm{d}Z_{u}$$
and $$f_g(r) = \int_0^t \phi'\left(  \int_0^{s-}
    g(s-v) \mathrm{d}Z_{v}  \right) \int_0^{s-}
    R(s-u,r) \mathrm{d}Z_{u} \mathrm{d}s.$$ 
Then the gradient of $l_t$ in $g$ is 
\begin{eqnarray*}
\nabla l_t(g) &=& f_g -  \sum_{i=1}^{N_t} \frac{ \phi'\left( \int_0^{\tau_i-} g(\tau_i-u)
    \mathrm{d}Z_{u}\right)}{\phi\left( \int_0^{\tau_i-} g(\tau_i-u)
    \mathrm{d}Z_{u}\right) } \eta_i.
\end{eqnarray*}
\end{prop}
 
The explicit derivation of the gradient above has several interesting consequences. First, a necessary 
condition for $g \in \Theta(D)^{\circ}$ to be a minimizer of the penalized
negative log-likelihood function is that $g$
solves $\nabla l_t(g) + 2 \lambda P g = 0$, 
which yields an integral equation in $g$. The integral equation is hardly solvable in any generality, but for
$\phi(x) = x + d$  it does provide the same information as Theorem \ref{thm:1} for interior 
minimizers -- that is, a minimizer must belong to the given finite dimensional subspace of $W^m_2$. 
The gradient can be used for descent algorithms. Inspired by the
gradient expression we propose a generic algorithm, Algorithm
\ref{alg:1}, for subspace approximations. We consider here only the
case where $D= \mathbb{R}$ so that $\Theta(D) = W^{m,2}$. 
The objective function that we attempt to minimize with Algorithm
\ref{alg:1} is 
$$\Lambda(g) = l_t(g) + \lambda ||P g||^2$$
with gradient $\nabla \Lambda(g) = \nabla l_t(g) + 2 \lambda P g.$ We
assume here that $\phi$ is continuously differentiable. To show a convergence result we need to introduce a
condition on the steps of the algorithm, and for this purpose we
introduce for $0 < c_1 < c_2 < 1$ and $\delta \in (0,1)$ fixed and $g \in W^{m,2}$ the
 subset  
{\small $$W(g) = \left\{ \tilde{g} \in W^{m,2} \backslash \{g\}
    \left|
    \ \begin{array}{rcl} \Lambda(\tilde{g}) - \Lambda(g)
& \leq &  c_1 \langle \nabla \Lambda(g), \tilde{g} - g \rangle  \\
 \langle \nabla \Lambda(\tilde{g}), \tilde{g} - g \rangle  & \geq & c_2 \langle  \nabla
 \Lambda(g), \tilde{g} - g \rangle   \\
- \langle  \nabla \Lambda (g), \tilde{g} - g\rangle  \ & \geq & \ \delta ||\nabla \Lambda (g)|| \,
||\tilde{g} - g|| \end{array}
 \right.\right\}$$} \\
The two first conditions determining $W(g)$ above are known as the \emph{Wolfe
  conditions} in the literature on numerical optimization, see
\citep{Nocedal:2006}. The third is an
\emph{angle condition}, which is automatically fulfilled if \mbox{$\tilde{g} -
g = - \alpha \nabla \Lambda(g)$} for $\alpha >
0$.  In Algorithm \ref{alg:1} we need to iteratively
choose $\hat{g}_h$, and we show that if $\nabla \Lambda(\hat{g}_{h-1})
\neq 0$ then under the
assumptions in Theorem \ref{thm:3} below  
\begin{equation} \label{eq:strongnonempty}
W(\hat{g}_{h-1}) \cap \text{span}\{\hat{g}_{h-1},\nabla
\Lambda(\hat{g}_{h-1})\} \neq
\emptyset,
\end{equation} 
which makes the iterative choices possible. 

\begin{figure}[!ht]
\begin{center}
\fbox{
\begin{minipage}[c]{.9\textwidth}
\begin{alg} \label{alg:1}
Initialize; fix $c_1,c_2$ with $0 < c_1 < c_2 < 1$ and $\delta \in (0,1)$, set
$$f_0(r) =  \int_0^t \int_u^t  R(s-u,r) \mathrm{d}s
\mathrm{d}Z_{u},$$
let $\hat{g}_0 \in \textnormal{span}\{\eta_1,\ldots,\eta_{N_t},f_0\}$ and set $h=1$.  
\begin{enumerate}
\item Stop if $\nabla \Lambda(\hat{g}_{h-1}) = 0$. Otherwise choose 
$$\hat{g}_{h} \in W(\hat{g}_{h-1}) \cap \textnormal{span}\{\eta_1,\ldots,\eta_{N_t},f_0,\ldots,f_{h-1}\}$$
where $W(g_{h-1})$ as defined above depends on $c_1$, $c_2$ and $\delta$.
\item Compute 
$$f_h(r) = \int_0^t \phi'\left(  \int_0^{s-}
    \hat{g}_{h}(s-v) \mathrm{d}Z_{v}  \right) 
   \int_0^{s-} R(s-u,r) \mathrm{d}Z_{u} \mathrm{d}s.$$
\item Set $h=h+1$ and return to 1. 
\end{enumerate}
\end{alg}
\end{minipage}
}
\end{center}
\end{figure}

\begin{thm} \label{thm:3} If $D=\mathbb{R}$, if $\phi$ is strictly positive, twice
  continuously differentiable and if the sublevel set 
$$\mathcal{S} = \{g \in \Theta(D) \mid \Lambda(g) \leq
\Lambda(\hat{g}_{0}) \}$$
is bounded then Algorithm \ref{alg:1} is globally convergent in the sense that 
$$||\nabla \Lambda(\hat{g}_h) || \rightarrow 0$$
for $h \to \infty$. 
\end{thm} 

If we, for instance, have strict convexity of $\Lambda$ then under
the assumptions in Theorem \ref{thm:3} we have a unique minimizer in
$\mathcal{S}$. Then we can strengthen the conclusion about
convergence and get weak convergence of $\hat{g}_h$ towards the
minimizer. In particular, we have the following corollary. 

\begin{cor} \label{cor:2} If there is a unique minimizer, $\hat{g}$, of $\Lambda$
  in $\mathcal{S}$ then under the assumptions in
  Theorem \ref{thm:3} 
$$\hat{g}_h(s) \rightarrow \hat{g}(s)$$
for $h \to \infty$ for all $s \in [0,t]$. 
\end{cor}

\subsection{Multivariate and additive models} \label{sec:additive}

We give in this section a brief treatment of how the setup in the
previous section extends to multivariate point processes and
to intensities given in terms of sums of linear filters. 

First we extend the models by considering additive intensities. We restrict the discussion to the situation
where $V = (W^{m,2})^d$ and $(Z_s)_{0 \leq s \leq t}$ is a
$d$-dimensional process. Perceiving $g \in V$ as a function
$g: [0,t] \to \mathbb{R}^d$ with coordinate functions in $W^{m,2}$ we
write 
$$\int_0^s g(s-u) \mathrm{d} Z_u = \sum_{j=1}^d \int_0^s g_j(s-u)
\mathrm{d} Z_{j,u}$$
and just as above, by Corollary
\ref{cor:fundcont}, 
$$g \mapsto X_sg := \int_0^s g(s-u) \mathrm{d} Z_u$$
is a continuous linear functional on $V$ when equipped with the product
topology. The inner product $\langle g,h\rangle  = \sum_{j=1}^d \langle g_j,h_j\rangle $ with
corresponding norm $||g||^2 = \sum_{j=1}^d ||g_j||^2$ turns
$V$ into a Hilbert space.

The negative log-likelihood function is given just as in the previous
section, but we will consider the more general penalization term  
$$J(g) = \lambda r(||Pg_1||^2,\ldots,||Pg_d||^2)$$
where $\lambda > 0$, $P$ is the orthogonal projection on
$\text{span}\{\gamma_1,\ldots,\gamma_l\}^{\perp}$ and $r:[0,\infty)^d
\to [0,\infty)$ is coordinate-wise increasing. Theorem \ref{thm:1}
easily generalizes with the following modification. If $\phi(x) = x+d$
then with
$$h_{i,j}(r) = \int_0^{\tau_i-} R^1(\tau_i-u,r) \mathrm{d} Z_{j,u}$$
for $i=1,\ldots,N_t$ and $j=1,\ldots,d$ a minimizer of the penalized
negative log-likelihood function has $j$'th
coordinate in the space spanned by $\gamma_1,\ldots,\gamma_l$, $h_{1,j},\ldots,h_{N_{t},j}$ and $f_j$ given by 
$$f_j(r) = \int_0^t \int_0^{s-} R^1(s-u,r)
\mathrm{d}Z_{j,u} \mathrm{d}s = \int_0^t \int_u^t 
R^1(s-u,r) \mathrm{d} s \mathrm{d}Z_{j,u}.$$
Proposition \ref{thm:2} also generalizes similarly and if $r$ is smooth,
for instance if $r(x_1,\ldots,x_d) = \sum_{j=1}^d x_j$, Algorithm
\ref{alg:1} generalizes as well. 

In the alternative, we can choose $r(x_1,\ldots,x_d) = \sum_{j=1}^d
\sqrt{x_j}$ leading to the penalty term 
$$J(g) = \lambda \sum_{j=1}^d ||Pg_j||,$$
which gives an infinite dimensional version of lasso. Since $r$ is not
differentiable, Algorithm \ref{alg:1} does not work directly. However,
a cyclical descent algorithm, as investigated thoroughly in
\citep{Hastie:2010} for the ordinary lasso, is implemented in \verb+ppstat+. Details
can be found in \citep{Hansen:2013b}.

The other direction of generalization is to the modeling of
multivariate point processes a.k.a. marked point processes with a discrete
mark space. The observation process is thus a multivariate counting
process $(N_{i,s})_{s \in [0,t]}$ for $i = 1, \ldots, p$ and we need to
specify separate intensities for each coordinate 
$$\lambda^i_s = \phi_i(X_s^i \beta_i)$$ 
for $\beta_i \in \Theta(D_i)$. With the coordinates being independent
homogeneous Poisson processes each with rate 1 under $P$, the negative log-likelihood becomes 
\begin{equation*}
  \sum_{i=1}^p \int_0^t \phi_i(X_s^i \beta_i) \mathrm{d}s
  -  \int_0^t  \log(\phi_i(X_s^i \beta_i)) \mathrm{d}N_{i,s}
\end{equation*}
for $\beta = (\beta_1, \ldots, \beta_p) \in \Theta(D_1) \times \ldots,
\times \Theta(D_p)$, see Theorem T.10, \citep{Bremaud:1981}. Since the $\beta_i$-parameters are variation
independent the negative log-likelihood is minimized by minimizing each term
separately. This caries over to the penalized negative log-likelihood if
the penalization function is of the form $J(\beta) = \sum_{i=1}^p
J_i(\beta_i)$, in which case the joint minimization reduces to $p$ separate
minimization problems -- one for each of the $p$ point processes. A typical
example is that $X_s = (N_{1,s}, \ldots, N_{p,s})$, that  
$$\beta_i = (g_{i1}, \ldots, g_{ip}) \in (W^{m,2})^p$$
and 
$$X_s^i \beta_i = \sum_{j=1}^p \int_0^{s-} g_{ij}(s-u) \mathrm{d}N_{j,u}.$$
Thus, the intensity for the $i$'th process has an additive
specification, as treated above, in terms of linear filters of the $p$
point processes.

\section{Example} \label{sec:examples}

In this section we work out some details for a more specific example
of Theorem \ref{thm:1}. For this we need a explicit choice of
inner product on $W^{m,2}$. Take
$$\mathcal{H}_1 = \{f \in  W^{m,2} \mid f(0) = Df(0) = \ldots = D^{m-1}f(0)=0\},$$
which we equip with the inner product 
$$\langle f,g\rangle  = \int_0^t D^mf(s) D^mg(s) \mathrm{d}s.$$ 
This turns $\mathcal{H}_1$ into a reproducing kernel Hilbert space
for $m \geq 1$ with reproducing kernel
 $R^1 : [0,t] \times [0,t] \to \mathbb{R}$ given as
$$R^1(s,r) = \int_0^{s \wedge r} \frac{(s-u)^{m-1}(r-u)^{m-1}}{((m-1)!)^2}
\mathrm{d} u,$$
see \citep{Wahba:1990}. Furthermore, define $\gamma_k(s) = s^{k-1}/(k-1)!$ for $k=1,\ldots,m$ and   
$$\mathcal{H}_0 = \textnormal{span}\{\gamma_1,\ldots,\gamma_m\},$$
which we equip with the inner product 
$$\langle \sum_{i} a_i \gamma_i, \sum_{j} b_j \gamma_j\rangle  = \sum_{i} a_ib_i,$$
so that $\gamma_1,\ldots,\gamma_m$ form an orthonormal basis for
$\mathcal{H}_0$. Then $\mathcal{H}_0$ is also a
reproducing kernel Hilbert space with reproducing kernel 
$R^0 : [0,t] \times [0,t] \to \mathbb{R}$ defined by
$$R^0(s,r) = \sum_{k=1}^m \gamma_k(s)\gamma_k(r).$$
Then the Sobolev space $W^{m,2} = \mathcal{H}_0 \oplus \mathcal{H}_1$ is a  
reproducing kernel Hilbert space with reproducing kernel 
$R(s,r) = R^0(s,r) + R^1(s,r)$, $\mathcal{H}_0 \perp \mathcal{H}_1$, 
and with $P$ the orthogonal projection onto $\mathcal{H}_1$, $PR = R^1$ and 
$$J(g) = \int_0^t (D^m g(s))^2 \mathrm{d}s.$$
It follows by the definition of $R$ that $R^1(s,\cdot)$ for fixed $s$ is
a piecewise polynomial of degree $2m-1$ with continuous derivatives of
order $2m-2$, that is, $R(s,\cdot)$ is an order $2m$ spline. We find 
that the $h_i$-functions for the basis in Theorem \ref{thm:1} are given
as stochastic integrals of order $2m$ splines.
% In the alternative, using Lemma \ref{lem:fubini}, we find that
% $$h_i(r) = \int_0^{\tau_i \wedge r} \int_0^{\tau_i - s} (\tau_i - s
% -u)^{m-1} \mathrm{d} Z_u  \frac{(r-s)^{m-1}}{((m-1)!)^2} \mathrm{d} s.$$
% Either way we conclude that explicit expressions for the basis 
% functions are obtainable if we can compute integrals of polynomials
% w.r.t. the semi-martingale. 

If $(Z_s)_{0 \leq s \leq t}$ itself is a counting process
  and $\phi(x) = x+d$ as in Theorem \ref{thm:1} we can give a more detailed
  description of the minimizer of (\ref{eq:minfunc}) over
  $\Theta(D)$. If $\sigma_1, \ldots, \sigma_{Z_t}$ denote the jump
  times for $(Z_s)_{0 \leq s \leq t}$ we find that 
$$h_i(r) = \sum_{j: \sigma_j < \tau_i} R^1(\tau_i - \sigma_j,r).$$
Collectively, the $h_i$ basis functions are order $2m$ splines with knots in 
$$\{ \tau_i - \sigma_j \mid i=1,\ldots,N_t, \ j: \sigma_j < \tau_i\}.$$
Due to (\ref{eq:fbasisinterchange}) the last basis function, $f$, is 
seen to be an order $2m+1$ spline with knots in 
$$\{ t- \sigma_j \mid i=1,\ldots,Z_t\}.$$
The cubic splines, $m=2$, are the splines mostly used in
practice. Here
$$R^1(s,r) = \int_0^{s \wedge r} (s-u)(r-u)
\mathrm{d} u = sr(s \wedge r) - \frac{(s+r)(s \wedge r)^2}{2} +
\frac{(s \wedge r)^3}{3}$$
and we can compute the integrated functions that enter in $f$ as
follows. If $t-u < r$
\begin{eqnarray*}
\int_u^t R^1(s-u,r) \mathrm{d} s = \int_0^{t-u} R^1(s,r) \mathrm{d} s 
%% &=& \frac{r(t-u)^3}{3} - \frac{r(t-u)^3}{6} - \frac{(t-u)^4}{8} + \frac{(t-u)^4}{12} \\
& = & \frac{r(t-u)^3}{6} - \frac{(t-u)^4}{24}
\end{eqnarray*} 
and if $t-u \geq s$
\begin{eqnarray*}
\int_u^t R^1(s-u,r) \mathrm{d} s &=& \int_0^{t-u} R^1(s,r) \mathrm{d} s =
\frac{3r^4}{24} + \int_r^{t-u} R^1(s,r) \mathrm{d} r \\
%% &=& \frac{3r^4}{24} + \left[\frac{s^2r^2}{4} - \frac{sr^3}{6} \right]^{t-u}_{r} \\
%% & = & \frac{3r^4}{24} + \frac{r^2(t-u)^2}{4} - \frac{r^3(t-u)}{6} -
%% \frac{r^4}{4} + \frac{r^4}{6} \\
& = & \frac{r^4}{24} + \frac{r^2(t-u)^2}{4} - \frac{r^3(t-u)}{6}. 
\end{eqnarray*}
Thus the function $f$ is a sum of functions, the $j$'th function being
a degree 4 polynomial on $[0,t-\sigma_j]$ and an affine function on
$(t-\sigma_j,t]$.

If $Z_s = N_s$ the process $(N_s)_{0\leq s \leq t}$ is under $Q_t$
known as a \emph{linear Hawkes process}, in which case the set of knots for
the $h_i$-functions is the collection of interdistances between the
points.

\section{Discussion}

The problem that initially motivated this paper was the estimation of the linear filter
functions entering in the specification of a non-linear Hawkes model
with an intensity specified as 
$$\phi\left( \sum_{j=1}^p \int_0^{s-} g_j(s-u) \mathrm{d}N_{j,u}\right)$$
where $N_j$ for $j=1,\ldots,p$ are counting processes, see
\citep{Bremaud:1996}. We have provided structural and algorithmic results for the
penalized maximum-likelihood estimator of $g_j$ in a Sobolev space, and
we have showed that these results can be established in a generality
where the stochastic integrals are with respect to any
c\`adl\`ag process. The representations of basis functions and the
gradient are useful for specific examples such as counting processes, but 
perhaps of limited analytic value for general processess. In practice we
can also only expect to observe a general process discretely and
numerical approximations to the integral representations and thus the
negative log-likelihood function must be used. If the process is 
coarsely observed it is unknown how reliable the resulting
approximation of the penalized negative log-likelihood function is. 

In this paper we relied on specific properties of Sobolev spaces to define
stochastic integrals pathwisely and to establish continuity of
certain linear functionals for general integrators. If we restrict
attention to pure jump process integrators the integrals are trivially
pathwisely defined and the required continuity properties follow by
elementary arguments for general reproducing kernel Hilbert
spaces with the minimal requirement that the reproducing kernel is 
continuous. See \cite{Hansen:2013b} for further details.

\section{Proofs} \label{sec:proofs}

\begin{proof} (Proposition \ref{prop:diff}).
If $\beta \in \Theta(D)^{\circ}$ and $\rho \in V$ then $\beta + q \rho
\in \Theta(D)^{\circ}$ for $q$ sufficiently small and
$$\partial_{q} \phi(X_s \beta + q X_s \rho) = \phi'(X_s
\beta +  q X_s \rho) X_s \rho.$$
Since $X_s \beta$ and $X_s \rho$ are bounded as functions of $s
\in [0,t]$ we have that 
$$(s, q) \mapsto X_s \beta + q X_s \rho$$
is bounded on $[0,t] \times [-\epsilon, \epsilon]$ and since $\phi'$
is assumed continuous we can bound $\partial_{q} \phi(X_s \beta + q
X_s \rho)$ uniformly by a constant on $[0,t] \times [-\epsilon,
\epsilon]$. This implies that we can interchange differentiation
w.r.t. $q$ and integration, thus 
\begin{eqnarray*}
\partial_{q} \int_0^t  \phi(X_s \beta + q
X_s \rho) \mathrm{d} s |_{q = 0} & =& \int_0^t  \partial_{q} \phi(X_s \beta + q
X_s \rho) |_{q = 0} \mathrm{d} s \\ 
&= & \int_0^t \phi'(X_s \beta) X_s \rho \mathrm{d} s.
\end{eqnarray*}
This gives the G\^ateaux derivative of the
first term in $\ell_t$.

For the second term note that $\ell_t(\beta) < \infty$ implies
that $\phi(X_{\tau_i-} \beta) > 0$, thus $\phi(X_{\tau_i-} \beta + q
X_{\tau_i-} \rho) > 0$, by continuity of $\phi$, for $q$ sufficiently
small, and
$$\partial_{q} \log \phi(X_{\tau_i-} \beta + q X_{\tau_i-} \rho) |_{q = 0} = \frac{\phi'(X_{\tau_i-}
\beta)}{\phi(X_{\tau_i-} \beta)} X_{\tau_i-} \rho.$$
This implies that the G\^ateaux derivative of the
second term in $\ell_t$ is
$$-\sum_{i=1}^{N_t} \frac{\phi'(X_{\tau_i-} \beta)}{\phi(X_{\tau_i-} \beta)}
X_{\tau_i-} = - \int_{0}^t \frac{\phi'(X_s \beta)}{\phi(X_s
  \beta)}  X_s \mathrm{d}N_s.$$

The second derivative is obtained similarly. 
\end{proof}

\begin{proof} (Proposition \ref{prop:diff2}).
The first derivative is found as above by decomposing the integration
interval $[0,t]$ into a finite number of closed intervals with the
roots of $X_s\beta$ as the end points.

For the second derivative note that the function 
$$H(s, q) = X_{s_1 + s}\beta + q X_{s_1 + s}\rho$$
is $C^1$ by assumption in $(-\delta, \delta) \times (-\epsilon, \epsilon)$ for
suitably small $\delta$ and $\epsilon$ and $s_1$ a root. Its
derivative w.r.t. $s$ in $(0, 0)$ is  
$$\partial_s H(0, 0) = \partial_s X_{s_1}\beta,$$
which is non-zero by assumption. Choosing $\delta$ small enough there
is, by the implicit function theorem, a $C^1$ function $s : (-\delta, \delta) \to (-\epsilon,
\epsilon)$ such that $H(s(q), q) = 0$. Assuming $\partial_s
X_{s_1}\beta >0$ and $X_{s_1}\rho \geq 0$ then $s(q) \geq 0$ for $q
\geq 0$ and  
$$1(H(s, q) > 0) - 1(H(s, 0) > 0) = - 1_{(0, s(q))}(s).$$
In particular, for an integrable function $h$ on $[-\delta,\delta]$ continuous in $0$
\begin{eqnarray*}
&& \frac{1}{q} \left(\int_{-\delta}^{\delta} h(s) 1(H(s, q) > 0) \mathrm{d} s -
\int_{-\delta}^{\delta} h(s) 1(H(s, 0) > 0) \mathrm{d} s\right) \\ 
&&  = - \frac{1}{q} \int_{0}^{s(q)} h(s) \mathrm{d} s \rightarrow
- s'(0) h(0)  =  \frac{X_{s_1}\rho}{\partial X_{s_1}\beta} h(0),
\end{eqnarray*}
from which the result follows.
% Is this complete?
\end{proof}

The Sobolev space $W^{m,2}$ was equipped with one particular inner
product denoted \mbox{$\langle \cdot , \cdot \rangle$} and
corresponding norm $||\cdot||$ in Section \ref{sec:examples}. An
alternative useful inner product on $W^{m,2}$ is
$$\langle f,g\rangle _{m} = \sum_{k=0}^m \int_0^t D^kf(s) D^kg(s) \mathrm{d}s$$
and the corresponding norm is given by
$$||f||_{m,2}^2 = \langle f,f\rangle _m = \sum_{k=0}^m \int_0^t D^kf(s)^2 \mathrm{d}s.$$
It is straight forward to show that $||\cdot||$ and $||\cdot||_{m,2}$
are equivalent norms, though the inner products give rise to different
reproducing kernels. For the theoretical arguments in this paper we will
use whichever norm is most convenient. Note that
the embedding $W^{m,2} \hookrightarrow W^{k,2}$ for $k < m$ is
continuous, which is straight forward using the norms
$||\cdot||_{m,2}$ and $||\cdot||_{k,2}$.  The continuity of the
embedding holds even when $k=0$ where $W^{0,2}=L_2$, which is not a
reproducing kernel Hilbert space.

We note that the characterizing property of a reproducing kernel Hilbert
space is that the function evaluations are continuous linear
functionals. If $\delta_s$ denotes the evaluation in $s$, that is,
$\delta_sf = f(s)$, then $R(s,\cdot)$ as a function in $W^{m,2}$
represents $\delta_s$ by 
$$f(s) = \langle f,R(s,\cdot)\rangle .$$
By Cauchy-Schwarz' inequality $||\delta_s|| = R(s,s)$ and since
$R$ is a continuous function of both variables, $R(s,s)$ is bounded
for $s$ in a compact set. 

We have already argued that stochastic integration of deterministic
functions from $W^{m,2}$ w.r.t. the a c\`adl\`ag process can be
defined by (\ref{eq:intebyparts}) as a pathwisely well defined linear
functional on $W^{m,2}$ for $m \geq 1$. We show that this functional
is continuous and subsequently that $X_s$ defined (\ref{eq:Xdef}) is
continuous.

\begin{lem} \label{lem:fundcont} Let $0 \leq s \leq t$. Then the
  linear functional $I_s : W^{1,2} \to \mathbb{R}$ defined by 
$$I_sh = \int_0^{s-} h(u) \mathrm{d}Z_u$$ 
is continuous. More precisely, we have the bound 
$$||I_s|| \leq |Z_{s-}|(1+s) + |Z_0| +  \left(\int_0^s Z_{u}^2\mathrm{d}u\right)^{1/2} < \infty.$$
\end{lem}

\begin{proof} Note that for $h \in W^{1,2}$ we have $$||h||^2 = |h(0)|^2 + ||h'||_2^2$$
and in particular 
$$||h'1_{[0,s]}||_{2}  \leq ||h'||_{2} \leq ||h||.$$ 

Using (\ref{eq:intebyparts}) and
Cauchy-Schwarz' inequality
\begin{eqnarray*}
|I_sh| & \leq & |h(s)Z_{s-}| + |h(0)Z_0| + \int_0^s |Z_uh'(u)|\mathrm{d} u \\
& \leq & |Z_{s-}|\,|h(s)| + |Z_0|\,|h(0)| + \left(\int_0^s Z_u^2\mathrm{d}u\right)^{1/2} ||h'1_{[0,s]}||_2 \\  
& \leq & \left(|Z_{s-}|\,||\delta_s|| + |Z_0|\,||\delta_0|| + \left(\int_0^s Z_u^2\mathrm{d}u\right)^{1/2}\right)||h|| \\
& \leq &  \left( |Z_{s-}|(1+s) + |Z_0| +  \left(\int_0^s Z_u^2\mathrm{d}u\right)^{1/2}\right)||h||,
\end{eqnarray*}
which shows the bound. Here we have used that for $m=1$ we
have $R(s,s) = 1+s$ and that $Z$ is c\`adl\`ag, hence bounded ($\omega$-wise)
and hence in $L_2([0,s])$ for any $s$. 
\end{proof}

In the following any function defined on $[0,t]$ is extended
to be 0 outside of $[0,t]$. Defining $\tau_s : W^{1,2} \to W^{1,2}$ by 
$$\tau_sg(u) = g(s) - \int_0^{u} g'(s - v) \mathrm{d} v
= \left\{ \begin{array}{ll}
g(s - u) & \quad \text{for } u \in [0,s] \\
g(0) & \quad \text{for } u \in (s,t] \end{array} \right.
$$
then $\tau_s$ is clearly linear and the linear functional $X_s$ defined by (\ref{eq:Xdef}) can
be expressed as $X_s = I_s \circ \tau_s$. 

\begin{lem} \label{lem:fundcont2} The linear operator $\tau_s :  W^{1,2} \to W^{1,2}$ is
  continuous with 
$$||\tau_s|| \leq \sqrt{(1+s)^2 + 1}.$$ 
\end{lem}

\begin{proof} We have that 
\begin{eqnarray*} 
||\tau_s(g)||^2 & = & |\tau_s(g)(0)|^2 + \int_0^t g'(s-v)^2 \mathrm{d}
v \\
& = & |g(s)|^2 + ||g'1_{[0, s]}||_2^2 \\
& \leq & ||\delta_s||^2||g||^2 + ||g||^2 = ((1 + s)^2 + 1) ||g||^2
\end{eqnarray*}
where we have used that $||\delta_s|| = R(s, s) = 1 + s$. Taking
square roots completes the proof.
\end{proof}

\begin{cor} \label{cor:fundcont} The linear functional $X_s : W^{m,2} \to \mathbb{R}$  
is continuous. Moreover, there is a real valued random
variable $C_{m,t}$ such that 
$$||X_s|| \leq C_{m,t}.$$
\end{cor}

\begin{proof} First consider $X_s = I_s \circ \tau_s : W^{1,2} \to
  \mathbb{R}$. Combining Lemma \ref{lem:fundcont} and Lemma
  \ref{lem:fundcont2} we find that
\begin{eqnarray*}
||X_s|| & \leq & ||I_s||||\tau_s|| \leq ||I_s|| \sqrt{(1 + s)^2 + 1} \\
& \leq & \sqrt{(1 + t)^2 + 1}  \left(\sup_{s \in [0,t]} |Z_s|(1+t) + |Z_0| + \left(\int_0^t
  Z_u^2\mathrm{d}u\right)^{1/2}\right) < \infty
\end{eqnarray*}
so the requested continuity of $X_s$ and the bound on $||X_s||$
follow. As the embedding
  $W^{m,2} \hookrightarrow W^{1,2}$ is continuous for $m
  \geq 1$ the $W^{m,2}$-norm of $X_s$ is
  bounded by a constant times the $W^{1,2}$-norm of $X_s$, and this
  completes the proof.
\end{proof}

We then turn to proving that $(X_s)_{0 \leq s \leq t}$ is continuous
from the left with right limits. 

\begin{lem} \label{lem:strongcont}
The map $s \mapsto \tau_s$ is strongly continuous
  from $[0,t]$ into the set of continuous linear operators on
  $W^{1,2}$, that is 
$$\lim_{\epsilon \to 0} ||\tau_s(g) - \tau_{s+\epsilon}(g)|| = 0$$
for all $g \in W^{1,2}$. 
\end{lem}

\begin{proof} Recall that even though $g'$ is initially defined on
  $[0,t]$, it is extended to be 0 outside of $[0,t]$, as mentioned above.  We then have that
\begin{eqnarray*}
||\tau_{s + \epsilon}(g) - \tau_{s}(g)||^2 & = & |g(s + \epsilon) -
g(s)|^2 + \int_0^t (g'(s + \epsilon - u) - g'(s - u))^2 \mathrm{d} u \\
& \leq & |g(s + \epsilon) - g(s)|^2 + ||g'(\cdot + \epsilon) - g'||_2^2
\end{eqnarray*}
with $||\cdot||_2$ denoting the $2$-norm on
$L_2(\mathbb{R})$ and with $g'(\cdot + \epsilon)$ denoting
the $\epsilon$-translation of $g'$. Since translation acts as a strongly
continuous (unitary) group on $L_2(\mathbb{R})$ we have that 
$$ ||g'(\cdot + \epsilon) - g'||_2^2 \to 0$$
for $\epsilon \to 0$ and continuity of $g$ ensures that also 
$|g(s + \epsilon) - g(s)|^2 \to 0$. For $s = 0$ or $s = t$ we only
consider limits from the right and from the left, respectively. 
\end{proof}

\begin{lem} \label{lem:intcaglad} The process $(I_s)_{0 \leq s \leq t}$ is continuous from
  the left with right limits in norm. 
\end{lem}

\begin{proof} Using (\ref{eq:intebyparts}) we have for $s \in (0,t]$ and
  $s \geq \epsilon > 0$ that
$$ |I_{s} h - I_{s-\epsilon} h| = \left|  
h(s)Z_{s-} - h(s-\epsilon) Z_{(s-\epsilon)-} + \int_{s-\epsilon }^{s} Z_u h'(u) \mathrm{d} u \right|,$$
and as in the proof of Lemma \ref{lem:fundcont} we get that 
$$ |I_s h - I_{s-\epsilon} h| \leq \left( ||Z_{s-} \delta_s - Z_{(s-\epsilon)-}
\delta_{s-\epsilon}||  +  \left(\int_{s-\epsilon}^{s}
Z_u^2 \mathrm{d} u\right)^{1/2}  \right) ||h||.$$
This shows that 
$$||I_{s} - I_{s-\epsilon}|| \leq ||Z_{s-} \delta_s - Z_{(s-\epsilon)-}
\delta_{s-\epsilon} ||  +  \left(\int_{s-\epsilon}^{s}
Z_u^2 \mathrm{d} u\right)^{1/2}$$
and letting $\epsilon \to 0$ the right hand side tends to 0 by an
application of dominated convergence and because 
$Z_{(s-\epsilon)-} \to Z_{s-}$ and  $\delta_{s-\epsilon} \to \delta_s$
in norm. 
% The norm convergence of $\delta_{s+\epsilon} \to \delta_s$ follows
% since the kernel is continuous and because $$||\delta_{s+\epsilon} -
% \delta_s|| = ||R(s+\epsilon,\cdot) - R(s,\cdot)|| =
% R(s+\epsilon,s+\epsilon) + R(s,s) - 2R(s+\epsilon,s)$$ 
This proves that the process is continuous from the left in norm. 

A similar argument shows that for $s \in [0,t)$ and $t-s \geq \epsilon > 0$
$$I_{s+\epsilon} \rightarrow I_s + (\Delta Z_s) \delta_s$$
for $\epsilon \to 0$ in norm where $\Delta Z_s = Z_{s} - Z_{s-}$. Thus
the process has limits from the right in norm.
\end{proof}

\begin{cor} \label{lem:caglad} The process $(X_s)_{0\leq s \leq t}$
  defined by (\ref{eq:Xdef}) is continous from the left and has limits
  from the right. 
\end{cor}

\begin{proof} We have for $g \in W^{1,2}$, $s \in (0,t]$ and $ s \geq \epsilon
  > 0$ that 
\begin{eqnarray*} 
|X_sg - X_{s - \epsilon} g| & = & |I_s(\tau_s(g)) -
I_{s-\epsilon}(\tau_{s-\epsilon}(g))| \\
& \leq &  |I_s(\tau_s(g)) -
I_s(\tau_{s-\epsilon}(g))| +  |I_s(\tau_{s-\epsilon}(g)) -
I_{s-\epsilon}(\tau_{s-\epsilon}(g))| \\ 
& \leq & ||I_s||||\tau_s(g) - \tau_{s-\epsilon}(g)|| + ||I_s -
I_{s-\epsilon}|| ||\tau_{s-\epsilon}(g)||.
\end{eqnarray*}
It follows from Lemma \ref{lem:strongcont} that $||\tau_s(g) - \tau_{s-\epsilon}(g)|| \to 0$
and from Lemma \ref{lem:intcaglad} that $||I_s - I_{s-\epsilon}|| \to
0$ for $\epsilon \to 0$.
Using the bounds in Lemma \ref{lem:fundcont2} and Corollary
\ref{cor:fundcont} on $||\tau_{s-\epsilon}||$ and $||I_s||$,
respectively, we conclude that the right hand side above converges to
0 for $\epsilon \to 0$ and thus that $X_{s-\epsilon}g \to X_sg$.  This
proves the left continuity. A similar argument shows that $X_s$ has
right limits.
\end{proof}

We may observe that 
$$\Delta X_s g = X_{s+} g - X_s g = \Delta Z_s g(0),$$
which shows that $s \mapsto X_s$ is (weak$^*$) continuous on the
subspace of $W^{1,2}$ where $g(0) = 0$. The map is in general continuous in
$s$ if and only if $Z_s$ is continuous in $s$.

We turn to the proof of Theorem \ref{thm:1} and for this purpose, as
well as for proving Proposition \ref{thm:2}, we will need the following
result.

\begin{lem} \label{lem:intfunc} Let $(H_t)_{t \geq 0}$ be a
 stochastic process with values in $V^*$, continuous from the left with
 right limits, and with $||H_s|| \leq C_t$ for
 $s \in [0, t]$ and $C_t$ a
 real valued random variable. Then the integral $\int_0^t H_{s} \mathrm{d}s$ defined by
\begin{equation} \label{eq:linfunc}
\beta \mapsto \int_0^t H_{s} \beta \mathrm{d}s
\end{equation}
is in $V^*$ with 
$$\left|\left|\int_0^t H_{s} \mathrm{d}s\right|\right| \leq t C_t.$$
\end{lem}

\begin{proof} The continuity requirements on $H_s$ implies that
  (\ref{eq:linfunc}) is well defined and clearly defines for a fixed $t \geq
  0$ a linear functional on $V$. Moreover, since $|H_{s} \beta| \leq
  ||H_{s}||\, ||\beta|| \leq C_t ||\beta||$
\begin{eqnarray*}
\left|\int_0^t H_{s}\beta \mathrm{d}s \right|  &\leq & 
\int_0^t \left| H_{s} \beta \right| \mathrm{d}s \\
& \leq & \int_0^t C_t \mathrm{d}s \
||\beta|| = t C_t ||\beta||.
\end{eqnarray*}
\end{proof}

\begin{proof} (Theorem \ref{thm:1}) When $\phi(x) = x+d$ we have that  

\begin{eqnarray*}
\ell_t(g)  &= & \int_0^t \int_0^{s-}
    g(s-u) \mathrm{d}Z_{u} + d \mathrm{d}s 
 -  \int_0^t \log\left(\int_0^{s-} g(s-u)
    \mathrm{d}Z_{u} + d \right) \mathrm{d}N_s \\
& = &  \int_0^t \int_0^{s-}
    g(s-u) \mathrm{d}Z_{u} \mathrm{d}s + td  -  \sum_{i=1}^{N_t} \log\left(\int_0^{\tau_i-} g(\tau_i-u)  \mathrm{d}Z_{u} + d\right).
\end{eqnarray*}
It follows from Corollary \ref{cor:fundcont} that
$$g \mapsto \int_0^{\tau_i-} g(\tau_i-u) \mathrm{d}Z_{u}$$
for $i=1,\ldots,N_t$ are continuous, linear functionals on $W^{m,2}$. The $i$'th of these continuous linear functionals is 
represented by $\eta_i \in W^{m,2}$ given as 
$$\eta_i(s) =  \int_0^{\tau_i-} R(\tau_i-u,s) \mathrm{d}Z_{u}$$
such that 
$$\langle \eta_i,g\rangle  =   \int_0^{\tau_i-} g(\tau_i-u) \mathrm{d}Z_{u}.$$ 
Hence $h_i = P\eta_i$. 

By combining Lemma \ref{lem:caglad} and Lemma \ref{lem:intfunc} we conclude that 
$$g \mapsto \int_0^t \int_0^{s-} g(s-u) \mathrm{d}Z_{u} \mathrm{d}s $$
is a continuous linear functional and $\eta$ is the representer given by 
$$\eta(r) =  \int_0^t \int_0^{s-}  R(s-u,r) \mathrm{d}Z_{u} \mathrm{d}s.$$ 
Hence  $f = P \eta$. 

Thus $\ell_t(g)$ is a function 
of a finite number of continuous, linear functionals on $W^{m,2}$,
$$\ell_t(g) = \langle \eta,g\rangle  - \sum_{i=1}^{N_t} \log (
\langle \eta_i,g\rangle) + td.$$
For $g \in \Theta(D) \subseteq W^{m,2}$, $g = g_0 + \rho$ with $\rho \in \textnormal{span}\{\gamma_1,\ldots,\gamma_l,h_1,\ldots,h_{N_t},f\}^{\perp}$, 
then  $\rho \perp \eta_i$ for $i=1,\ldots,N_t$, $\rho \perp \eta$, $P\rho=\rho$ and 
\begin{eqnarray*}
\ell_t(g) + \lambda ||Pg||^2 & = & \langle \eta,g\rangle  - \sum_{i=1}^{N_t}
\log (\langle \eta_i,g\rangle)  + td + \lambda ||Pg||^2 \\
& = & \langle \eta,g_0\rangle  - \sum_{i=1}^{N_t} \log (\langle \eta_i,g_0\rangle)
+ td + \lambda ||Pg_0||^2 + \lambda ||\rho||^2 \\
& \geq & \ell_t(g_0) + \lambda ||Pg_0||^2  \\
\end{eqnarray*}
with equality if and only if $\rho = 0$. Thus a minimizer of $\ell_t(g) + \lambda ||Pg||^2$ over $\Theta(D)$ must 
be in $\textnormal{span}\{\gamma_1,\ldots,\gamma_l,h_1,\ldots,h_{N_t},f\}$.
\end{proof}

We have used the Fubini
theorem below to give an
alternative representation of the
basis function $f$ from Theorem \ref{thm:1}. The result is a
consequence of Theorem 45 in \citep{Protter:2005} when the integrator
is a semi-martingale. With the pathwise 
definition of stochastic integrals, as given by
(\ref{eq:intebyparts}), we give an elementary proof. 

\begin{lem} \label{lem:fubini} With $(Z_s)_{0 \leq s \leq t}$ a
  c\`adl\`ag process, $(Y_s)_{0 \leq s \leq t}$ a c\`agl\`ad
  process and $g \in W^{1,2}$ then
$$ \int_0^t Y_s  \int_0^{s-}
    g(s-u) \mathrm{d}Z_{u} \mathrm{d}s =\int_0^t \int_u^t Y_s  g(s-u) \mathrm{d}s \mathrm{d}Z_{u}.$$
\end{lem}

\begin{proof}
Using (\ref{eq:intebyparts}) and Fubini
\begin{eqnarray*} 
 \int_0^t Y_s  \int_0^{s-}
    g(s-u) \mathrm{d}Z_{u} \mathrm{d}s & = &  g(0) \int_0^t Z_{s} Y_s
    \mathrm{d}s - Z_0 \int_0^t g(s)Y_s \mathrm{d}s \\ && \hspace{2cm}
+ \int_0^t Y_s \int_0^{s} Z_u g'(s-u)\mathrm{d}u \mathrm{d} s \\
& = &  g(0) \int_0^t Z_{s} Y_s \mathrm{d}s - Z_0 \int_0^t g(s)Y_s \mathrm{d}s \\ && \hspace{2cm}
+ \int_0^t Z_u \int_u^t Y_s g'(s-u)\mathrm{d}s \mathrm{d} u.
\end{eqnarray*}

To use (\ref{eq:intebyparts}) on the right hand side above we 
need to verify that the integrand is sufficiently regular. Defining  
$$G(u) =  \int_u^t Y_s  g(s-u) \mathrm{d}s$$
for $g \in W^{1,2}$ then $G$ is weakly differentiable with derivative
$$G'(u) = - \int_u^t Y_s  g'(s-u) \mathrm{d}s - Y_u g(0),$$
which is verified simply by checking that $G(u) = - \int_u^t G'(v)
\mathrm{d}v$. Using this, we get for the right hand side above that
\begin{eqnarray*}
\int_0^t \underbrace{\int_u^t Y_s  g(s-u) \mathrm{d}s}_{G(u)} \mathrm{d}Z_{u} & = & G(t)Z_{t} - G(0)Z_0 - \int_0^t Z_u G'(u) \mathrm{d}u \\
& = & - G(0)Z_0 +  \int_0^t Z_u  \left[\int_u^t Y_s  g'(s-u) \mathrm{d}s + Y_u g(0) \right]\mathrm{d} u \\ 
& = & g(0) \int_0^t Z_{s} Y_s \mathrm{d}s -  Z_0 \int_0^t g(s)Y_s
\mathrm{d}s \\ && \hspace{2cm} + \int_0^t Z_u  \int_u^t Y_s  g'(s-u) \mathrm{d}s \mathrm{d} u.
\end{eqnarray*}
\end{proof}

\begin{proof} (Proposition \ref{thm:2}) The G\^ateaux derivative of $l_t$ in the
  direction of $h \in W^{m,2}$ for $g \in \Theta(D)^{\circ}$ is by Proposition \ref{prop:diff}
\begin{eqnarray*}
D l_t(g)h &=& \int_0^{t} \phi'\left(  \int_0^{s-}
    g(s-u) \mathrm{d}Z_{u}  \right) \int_0^{s-} h(s-u) \mathrm{d}Z_{u}  \mathrm{d}s \\ 
&& -   \int_0^t \frac{ \phi'\left( \int_0^{s-} g(s-u)
    \mathrm{d}Z_{u}\right)}{\phi\left( \int_0^{s-} g(s-u)
    \mathrm{d}Z_{u}\right) }\int_0^{s-} h(s-u)
    \mathrm{d}Z_{u} \mathrm{d}N_s.
\end{eqnarray*}
Now just as in the proof of Theorem \ref{thm:1}, using Lemma
\ref{lem:caglad} and Lemma \ref{lem:intfunc} with 
$$H_sh  =  \phi'\left(  \int_0^{s-}
    g(s-u) \mathrm{d}Z_{u}  \right) \int_0^{s-} h(s-u)
  \mathrm{d}Z_{u},$$ 
the first term is a continuous linear functional on $W^{m,2}$
with representer $f_g$.
Moreover, with $\eta_i$ as defined in Proposition \ref{thm:2} the second term above is seen to be a
continuous linear functional on  $W^{m,2}$ with representer 
$$\zeta_g = \sum_{i=1}^{N_t} \frac{ \phi'\left( \int_0^{\tau_i-} g(\tau_i-u)
    \mathrm{d}Z_{u}\right)}{\phi\left( \int_0^{\tau_i-} g(\tau_i-u)
    \mathrm{d}Z_{u}\right) } \eta_i.$$
In conclusion, the gradient of $l_t$ in $g$ is $\nabla l_t(g) = f_g - \zeta_g.$
\end{proof}

\begin{lem} If $D=\mathbb{R}$ and $\phi$ is strictly positive, twice
  continuously differentiable then 
the gradient $\nabla \Lambda : W^{m,2} \to W^{m,2}$ is
Lipschitz continuous on any bounded set. 
\end{lem}

\begin{proof} Let $B(0,L)$ denote the ball with radius $L$ in
  $W^{m,2}$. Corollary \ref{cor:fundcont} shows that $|X_sg|
  \leq C_{m,t} ||g||$. This means that there is an $M > 0$ such that $X_sg \in
  [-M,M]$ for all $g \in B(0,L)$ and $s \in [0,t]$. Since $\phi$ is twice continuously differentiable we
  have that $\phi'$ is Lipschitz continuous on $[-M,M]$ with Lipschitz constant $K$, say. With $f_g$ for $g \in W^{m,2}$ as in
  Theorem \ref{thm:3} we find that for $g, g' \in W^{m,2}$
\begin{eqnarray*}
f_g - f_{g'} & = & \int_0^t \phi'( X_sg) -  \phi'( X_sg')
\int_0^{s-} R^1(s-u,\cdot)\mathrm{d} Z_u  \mathrm{d} s \\
\end{eqnarray*}
and as above, by the isometric isomorphism that identifies $W^{m,2}$ with
its dual, we get by Lemma \ref{lem:intfunc} that if also $g, g' \in
B(0,L)$ then 
\begin{eqnarray*}
||f_g - f_{g'}|| & \leq &  C_{m,t} \int_0^t | \phi'( X_sg)
- \phi'( X_sg')| \mathrm{d} s \\
& \leq & \underbrace{K t C_{m,t}^2}_{C_1} ||g - g'||.
\end{eqnarray*}
Since $\phi$ is strictly positive and twice continuously
differentiable, the function $x \mapsto \phi'(x)/\phi(x)$ is Lipschitz
continuous on $[-M,M]$ with Lipschitz constant $K'$, say. Then for
$g,g' \in B(0,L)$
\begin{eqnarray*}
\left|\left|\sum_{i=1}^{N_t} \frac{ \phi'\left( X_{\tau_i} g
  \right)}{\phi\left( X_{\tau_i} g \right) } \eta_i - 
\sum_{i=1}^{N_t} \frac{ \phi'\left( X_{\tau_i} g'
  \right)}{\phi\left( X_{\tau_i} g' \right) } \eta_i\right|\right|
& \leq &\sum_{i=1}^{N_t} \left|\frac{ \phi'\left( X_{\tau_i} g
  \right)}{\phi\left( X_{\tau_i} g \right) } - 
 \frac{ \phi'\left( X_{\tau_i} g'
  \right)}{\phi\left( X_{\tau_i} g' \right) } \right| \,
||\eta_i|| \\
& \leq &K' \sum_{i=1}^{N_t} ||X_{\tau_i}||\, ||g-g'||\, ||\eta_i|| \\
& \leq &\underbrace{K' \left(\sum_{i=1}^{N_t} ||X_{\tau_i}||\, ||\eta_i||\right)}_{C_2} ||g-g'||.
\end{eqnarray*}
By Proposition \ref{thm:2} we have showed that the gradient $\nabla l_t$ is Lipschitz
continuous on the bounded set $B(0,L)$ with Lipschitz
constant $C=C_1+C_2$. Since $\nabla \Lambda  = \nabla l_t + 2 \lambda
P$ and $2 \lambda P$ is linear this proves that $\nabla \Lambda$ is Lipschitz
continuous on bounded sets. 
\end{proof}

\begin{proof} (Theorem \ref{thm:3}) We prove first by induction that it is possible
  to iteratively choose $\hat{g}_h$ as prescribed in Algorithm
  \ref{alg:1}. The induction start is given by assumption. 

Assume that $\hat{g}_h$ is chosen as in Algorithm \ref{alg:1}. Since
$\Lambda : W^{m,2} \to \mathbb{R}$ is continuous and 
$$\mathcal{S}_h := \{g \in W^{m,2} \mid \Lambda(g) \leq 
\Lambda(\hat{g}_{h}) \} \subseteq \mathcal{S}$$
is bounded by assumption we find that $\Lambda$ is bounded below along
the ray $\hat{g}_h - \alpha \nabla \Lambda(\hat{g}_h)$ for $\alpha >
0$. If $\nabla \Lambda(\hat{g}_h) \neq 0$ we can proceed exactly as in
the proof of Lemma 3.1 in
\citep{Nocedal:2006}, and there exists $\alpha > 0$ such that 
$$\tilde{g}_{h+1} = \hat{g}_h - \alpha \nabla \Lambda(\hat{g}_h) \in \mathcal{S}_h$$ 
fulfills the two Wolfe conditions:
\begin{eqnarray*}
\Lambda(\tilde{g}_{h+1}) & \leq  & \Lambda(\hat{g}_h) - c_1 \alpha ||\nabla \Lambda(\hat{g}_h)||^2 \\
 \langle \nabla \Lambda(\tilde{g}_{h+1}), \nabla \Lambda(\hat{g}_h)\rangle  & \leq
&  c_2 ||\nabla \Lambda(\hat{g}_h)||^2.
\end{eqnarray*}
Since $\tilde{g}_{h+1} - \hat{g}_h = - \alpha \nabla
\Lambda(\hat{g}_h) \neq 0$,  and since $\hat{g}_h \in \text{span}\{\eta_1,\ldots,\eta_{N_t},f_0,\ldots,f_{h-1}\}$ and $\nabla \Lambda(\hat{g}_h)
\in  \text{span}\{\eta_1,\ldots,\eta_{N_t},f_0,\ldots,f_{h}\}$ we find that 
$$\tilde{g}_{h+1} \in  W(\hat{g}_h) \cap
\text{span}\{\eta_1,\ldots,\eta_{N_t},f_0,\ldots,f_{h}\}$$
and the set on the right hand side is, in particular,
non-empty. This proves that it is possible to
iteratively choose $\hat{g}_h$ as in Algorithm \ref{alg:1}. 

For the entire sequence $(\hat{g}_h)_{h \geq 0}$ we get from the second Wolfe
condition together with the Cauchy-Schwarz inequality and Lipschitz
continuity of $\nabla \Lambda$ on $\mathcal{S}$ that 
\begin{eqnarray*}
(c_2 - 1) \langle  \nabla \Lambda(\hat{g}_h), \hat{g}_{h+1} -
\hat{g}_h\rangle  & \leq & \langle \nabla \Lambda(\hat{g}_{h+1}) - \nabla
\Lambda(\hat{g}_{h}), \hat{g}_{h+1} - \hat{g}_h\rangle  \\
& \leq &  C ||\hat{g}_{h+1} - \hat{g}_h||^2,
\end{eqnarray*}
which implies that 
$$||\hat{g}_{h+1} - \hat{g}_h|| \geq \frac{(c_2-1)}{C}\frac{\langle  \nabla
\Lambda(\hat{g}_h), \hat{g}_{h+1} - \hat{g}_h\rangle }{||\hat{g}_{h+1}
- \hat{g}_h||}.$$ 
Note that, by the angle condition, the inner product above is
\emph{strictly negative} when $\nabla \Lambda(\hat{g}_n) \neq 0$, and since $c_2 < 1$ this lower bound is actually
always non-trivial.
Combining the angle condition with the first Wolfe
condition gives that  
\begin{eqnarray*}
\Lambda(\hat{g}_{h+1}) &\leq & \Lambda(\hat{g}_h) + c_1
||\hat{g}_{h+1} - \hat{g}_h|| \frac{\langle  \nabla
\Lambda(\hat{g}_h), \hat{g}_{h+1} - \hat{g}_h \rangle }{||\hat{g}_{h+1}
- \hat{g}_h||} \\
& \leq & \Lambda(\hat{g}_h) - \frac{c_1(1-c_2)}{C} \frac{\langle  \nabla
\Lambda(\hat{g}_h), \hat{g}_{h+1} - \hat{g}_h\rangle ^2 }{||\nabla
\Lambda(\hat{g}_h)||^2 ||\hat{g}_{h+1} - \hat{g}_h||^2} ||\nabla
\Lambda(\hat{g}_h)||^2 \\
& \leq & \Lambda(\hat{g}_h) - \frac{c_1(1-c_2)\delta^2}{C} ||\nabla \Lambda(\hat{g}_h)||^2.
\end{eqnarray*}
By induction
$$\Lambda(\hat{g}_{h+1}) \leq \Lambda(\hat{g}_{0}) -
\frac{c_1(1-c_2)\delta^2}{C} \sum_{k=0}^h ||\nabla
\Lambda(\hat{g}_k)||^2.$$
To finish the proof we need to show that $\Lambda$ is bounded below on
$\mathcal{S}$, because then the inequality above implies that  
$$||\nabla \Lambda(\hat{g}_h)|| \rightarrow 0$$
for $h \to \infty$. To show that $\Lambda$ is bounded below we observe
that 
\begin{eqnarray*}
\Lambda(g) & \geq&  - \int_0^t \log(\phi\left(\int_0^{s-} g(s-u)
  \mathrm{d} Z_u\right) )\mathrm{d}N_s \\
& = & - \sum_{i=1}^{N_t} \log( \phi\left(\int_0^{s-} g(\tau_i-u)
  \mathrm{d} Z_u\right) ) \\
& = & - \sum_{i=1}^{N_t} \log( \phi(\langle \eta_i,g\rangle )).
\end{eqnarray*}
Since this lower bound as a function of $g$ is weakly continuous, and
since a bounded set is weakly compact by reflexivity of a Hilbert
space and Banach-Alaoglu's theorem, 
% Theorem 2.5.2, \citep{GKP:1989}, 
we have proved that $\Lambda$ is bounded below 
on the bounded set $\mathcal{S}$.
\end{proof}

For the proof of Corollary \ref{cor:2} we need the following lemma.

% This may not be correct, and as it turns out, it is not needed.
% \begin{lem} If $g_n \overset{w}{\rightarrow} g$ in $W^{1,2}$
%   then $\int_0^t |g_n'(s)| \mathrm{d} s$ is a bounded sequence. 
% \end{lem}

% \begin{proof} First note that with $\gamma_2(t) = t$ then for $f \in W^{1,2}$
% $$\langle \gamma_2, f \rangle = \gamma_2(0) f(0) + \int_0^t
% \gamma_2'(s) f'(s) \mathrm{d} s = \int_0^t f'(s) \mathrm{d} s.$$
% If $f'_+$ and $f'_-$ denote the positive and negative part of $f'$,
% respectively, then $f'_+, f'_- \in L_2$ and with 
% $$f_+(s) = \int_0^s f'_+(u) \mathrm{d} u \quad \text{and} \quad 
% f_-(s) = \int_0^s f'_-(u) \mathrm{d} u$$
% we see that $f_+, f_- \in W^{1,2}$. Since $|f'(s)| = f'_+(s) + f'_-(s)$
% we get that 
% $$\int_0^t |f'(s)| \mathrm{d} s = \int_0^t f'_+(s) \mathrm{d} s +
% \int_0^t f_-'(s) \mathrm{d} s = \langle \gamma_2, f_+ \rangle +
% \langle \gamma_2, f_- \rangle = \langle \gamma_2, f_+ + f_-\rangle.$$
% Thus, the $1$-norm of $f'$ is given as a linear functional of $f_+ +
% f_-$. Here is a problem, because even if $g_n$ converges weakly the
% functions $(g_n)_+$ and $(g_n)_-$ do not automatically converge ... 
% \end{proof}

\begin{lem} \label{lem:weakweakcont} If $\phi$ is strictly positive
  and continuously differentiable the map $g \mapsto \nabla
  \Lambda(g)$ is sequentially weak-weak continuous.
\end{lem}

\begin{proof} By definition of the weak topology we need to show that 
$$g \mapsto \langle  \nabla \Lambda(g), h\rangle  = \langle \nabla l_t(g), h\rangle  + 2\lambda \langle 
P g, h\rangle  = Dl_t(g)h + 2\lambda \langle Pg, h \rangle$$
is weakly continuous for all $h \in W^{1,2}$. Clearly $g
\mapsto \langle  P g, h\rangle =\langle  g, P h\rangle $ is weakly continuous so we can restrict
our attention to $g \mapsto D l_t(g)h$. We use
Proposition \ref{prop:diff}, and observe that the continuous linear functional
$$g \mapsto X_sg = \int_0^{s-} g(s-u) \mathrm{d}Z_u$$
for fixed $s$ is weakly continuous by the definition of the weak
topology. We conclude
directly from this that
$$ g \mapsto \sum_{i=1}^{N_t} \frac{ \phi'\left( X_{\tau_i} g
    \right)}{\phi\left( X_{\tau_i} g \right) } X_{\tau_i}h$$
is weakly continuous as $\phi$ is assumed strictly positive and
continuously differentiable. To handle the second term in the
derivative assume that $g_n \overset{w}{\rightarrow} g$ for $n \to \infty$,
in which case 
$$X_s g_n \rightarrow X_s g$$
for all $s \in [0,t]$. By the uniform boundedness principle
(the Banach-Steinhaus theorem) the weakly convergent sequence $(g_n)_{n
  \geq 1}$ is bounded in $W^{m,2}$. Then it follows
from the bound on $||X_s||$ in Corollary \ref{cor:fundcont} that
$$\sup_n \sup_{s \in [0,t]} |X_s g_n| \leq C_{m,t}
\sup_{n} ||g_n|| < \infty.$$

Since $\phi'$ is continuous the pointwise convergence of 
$$\phi'\left(X_s g_n \right) X_sh \to
\phi'\left(X_s g \right) X_sh$$
for $s \in [0,t]$ is dominated by a constant, which is integrable over $[0,t]$. Hence 
\begin{eqnarray*}
\int_0^t  \phi'\left(X_s g_n \right) X_sh \mathrm{d} s 
& \rightarrow & \int_0^t \phi'\left(X_s g \right) X_sh \mathrm{d} s
\end{eqnarray*}
for $n \to \infty$. 
\end{proof}

Whether $\nabla \Lambda$ is actually weak-weak continuous on
$W^{m,2}$ and not just sequentially weak-weak continuous is
not of our concern. Since bounded sets in the Hilbert space
$W^{m,2}$ are metrizable in the weak topology, $\nabla \Lambda$
is weak-weak continuous on every bounded set. In the following proof 
weak-weak continuity on a bounded set suffices. 

\begin{proof} (Corollary \ref{cor:2}) By assumption, $\hat{g} \in \mathcal{S}$ is the
  unique solution to \mbox{$\nabla \Lambda (g) = 0$.}
The bounded set $\mathcal{S}$ is weakly compact as argued above and
the weak topology is, moreover, metrizable on $\mathcal{S}$ since
$W^{m,2}$ is separable. Therefore any subsequence of $(\hat{g}_h)_{h \geq 0}$ has a subsequence that
converges weakly in $\mathcal{S}$, necessarily towards a limit with
vanishing gradient by Lemma \ref{lem:weakweakcont}. Uniqueness of $\hat{g}$ implies that
$(\hat{g}_h)_{h \geq 0}$ itself is weakly convergent with limit
$\hat{g}$. The proof is completed by noting that weak convergence in a reproducing 
kernel Hilbert space implies pointwise convergence. 
\end{proof}

\bibliographystyle{model1b-num-names}
\bibliography{../texbib-1}

\end{document}